\documentclass[12pt]{article}
\usepackage{amsfonts}
\usepackage{amsmath}
\usepackage{amssymb}
\usepackage[english]{babel}
\usepackage[T1]{fontenc}
\usepackage{graphicx}
\usepackage{subfigure}
\newtheorem{prop}{Proposition}[section]
\newtheorem{rem}[prop]{Remark}

\newtheorem{theo}[prop]{Theorem}
\newtheorem{corr}[prop]{Corollary}

\usepackage{epsfig}
\usepackage{latexsym}
\usepackage{amstext}
\usepackage{epsfig, graphicx}
\usepackage{epsfig}
\usepackage{latexsym}
\usepackage{amstext}
\usepackage{amssymb}
\usepackage{graphics}
\usepackage{enumerate}

\newcommand{\beq}{\begin{eqnarray}}
\newcommand{\beqq}{\begin{eqnarray*}}
\newcommand{\eeq}{\end{eqnarray}}
\newcommand{\eeqq}{\end{eqnarray*}}

\def\QED{\quad\hbox{\hskip 4pt\vrule width 5pt height 6pt depth 1.5pt}}

\title{The Mean First Rotation Time of a planar polymer}
\author{ S. Vakeroudis\thanks{Laboratoire de Probabilit\'{e}s et Mod\`{e}les
Al\'{e}atoires (LPMA) CNRS : UMR7599,  Universit\'{e} Pierre et Marie
Curie - Paris VI. } \thanks{D\'epartement de Math\'ematiques et de
Biologie, Ecole Normale Sup\'erieure, 46 rue d'Ulm 75005 Paris,
France. D. H. is supported by HFSP and ERC-SG.} \
M. Yor$^{\ast}$\thanks{Institut Universitaire de
France.} \
D. Holcman$^{\dag}$\thanks{Department of Applied Mathematics, Tel-Aviv University, Tel-Aviv 69978,
Israel.} }
\date{\today}
\begin{document}

\maketitle
\begin{abstract}
We estimate the mean first time, called the mean rotation time
(MRT), for a planar random polymer to wind around a point. This polymer is
modeled as a collection of $n$ rods, each of them being parameterized by a
Brownian angle. We are led to study the sum of i.i.d. imaginary exponentials
with one dimensional Brownian motions as arguments. We find that
the free end of the polymer satisfies a novel stochastic equation
with a nonlinear time function. Finally, we obtain  an asymptotic formula for the MRT, whose leading order term
depends on $\sqrt{n}$ and, interestingly, depends weakly on the mean initial
configuration. Our analytical results are confirmed by Brownian
simulations.
\end{abstract}

\section{Introduction}
\renewcommand{\thefootnote}{\fnsymbol{footnote}}
This paper focusses mainly on some properties of a planar
polymer motion and in particular, on the mean time that a rotation is
completed around a fixed point. This mean rotation time (MRT)
provides a quantification for the transition between a free two
dimensional Brownian motion and a restricted motion. To
study the distribution of this time, we shall use a simplified model of a
polymer made of a collection of $n$ two dimensional connected rods
with Brownian random angles (Figure \ref{figure1}). To estimate the
mean rotation time, we shall fix one end of the polymer. We
shall examine how the MRT depends on various parameters such as
the diffusion constant, the number of rods or their common length. Using
some approximations and numerical simulations, the mean
time for the two polymer ends to meet was estimated in dimension
three \cite{WiF74,PZS96}.

Although the windings of a (planar) Brownian motion, that is the
number of rotations around one point in dimension 2 or a line in
dimension three and its asymptotic behavior have been
studied quite extensively \cite{Spi58,PiY86,LeGY87,LeG90,ReY99}, little seems to be
known about the mean time for a rotation to be completed for the
first time. For an Ornstein-Uhlenbeck process, the first time that it hits the boundary of a
given cone was recently estimated \cite{Vak10}.

The paper is organized as follows: in section 2, we present the
polymer model. In section 3, we study a sum of i.i.d.
exponentials of one dimensional Brownian motions. Interestingly,
using a Central Limit Theorem, we obtain a new stochastic equation
for the limit process. This equation describes the motion of the
free polymer end. In section 4, we obtain our main result which
is an asymptotic formula for $E\left[\tau_{n}\right]$ the MRT  when
the polymer is made of $n$ rods of equal lengths $l_{0}$ with the
first end fixed at a distance $L$ from the origin and the Brownian
motion is characterized by its rotational
diffusion constant $D$. We find that the MRT depends logarithmically  on the mean initial configuration and
for $nl_{0} >> L$ and $n\geq 3$, the leading order term is given
by:
\begin{enumerate}
    \item { for a general initial configuration:
\beqq
E\left[\tau_{n}\right] \approx  \frac{\sqrt{n}}{8D}
\left[2\ln \left( \frac{E \left(\sum^{n}_{k=1} e^{\frac{i}{\sqrt{2D}}
\theta_{k}(0)}\right)}{\sqrt{n}}\right)+Q\right],
\eeqq
where $\left(\theta_{k}(0),1\leq k\leq n\right)$ are the
initial angles of the polymer and $Q\approx9.56$, }
\item { for an initially  stretched polymer:
\beqq
E\left[\tau_{n}\right] \approx  \frac{\sqrt{n}}{8D} \; \left(\ln (n)+Q\right),
\eeqq } \\
    \item {for an average over uniformly distributed initial angles:
\beqq
E\left[\tau_{n}\right] \approx  \frac{\sqrt{n}}{8D} \; \tilde{Q},
\eeqq
where $\tilde{Q}\approx9.62$.
}
\end{enumerate}
We confirm our analytical results with generic Brownian
simulations with normalized parameters . Finally, in section 5 we discuss some related
open questions.

\section{Stochastic modeling of a planar polymer}
\begin{figure}
\centering
\begin{tabular}{cc}
\includegraphics[width=0.55\textwidth]{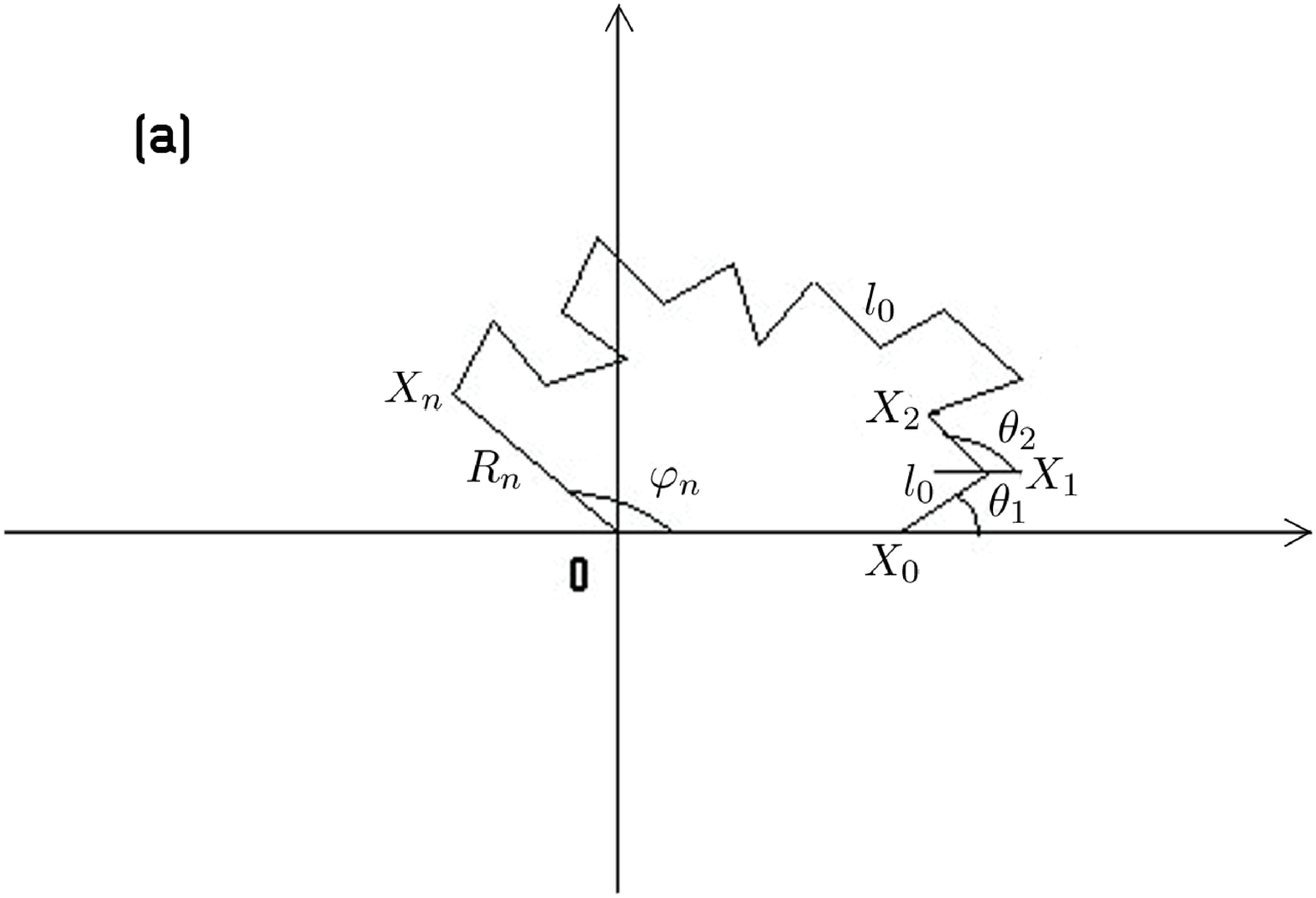} &
\includegraphics[width=0.55\textwidth]{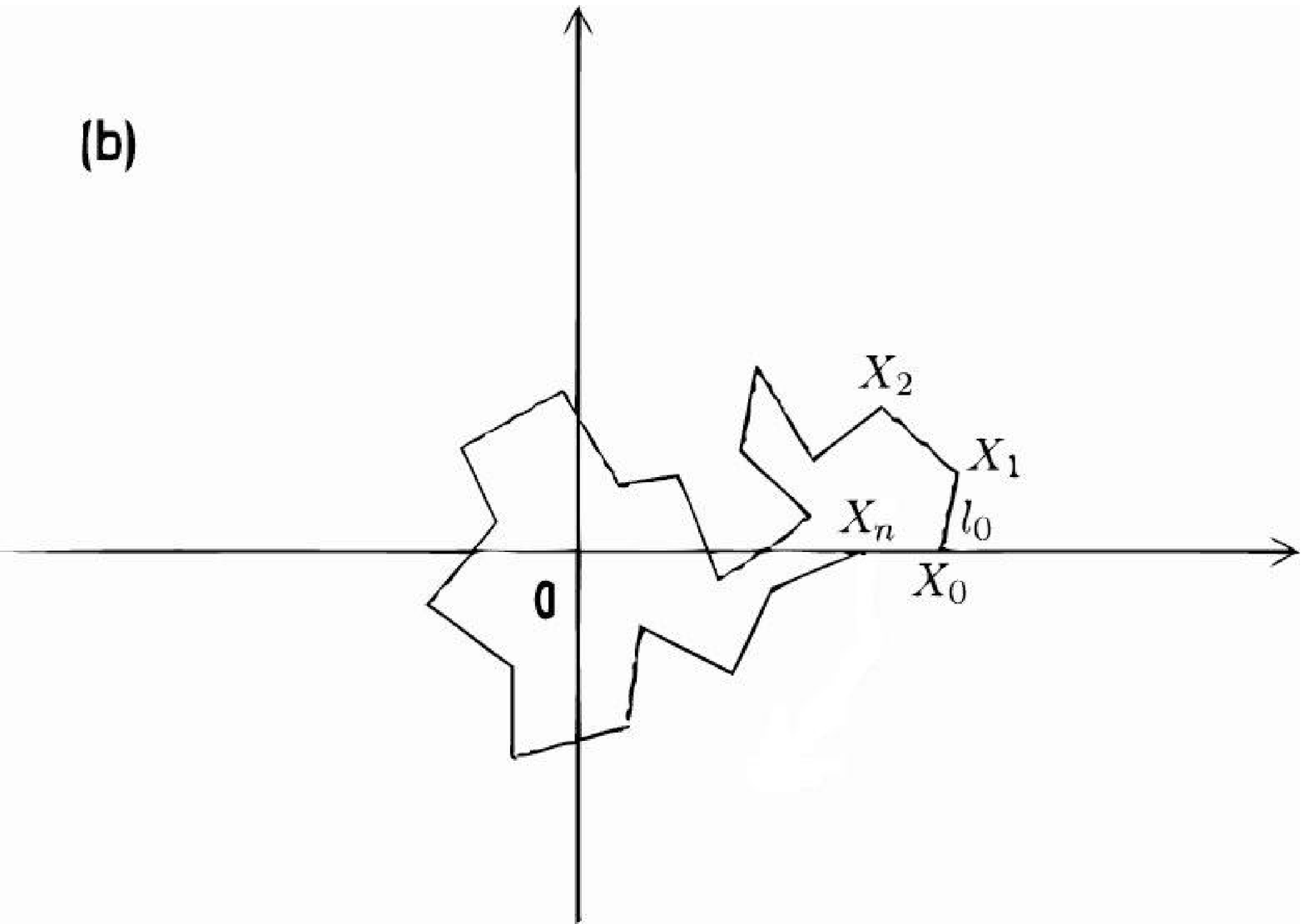} \\
\end{tabular}
    \caption{ {\bf Schematic representation of a planar polymer winding
     around the origin}.(a) A random configuration,
        (b) MRT when the $n$-th bead reaches $2\pi$.}
   \label{figure1}
\end{figure}

Various models are available to study polymers: the Rouse model
consists of a collection of beads connected by springs, while more
sophisticated models account for bending, torsion and specific
mechanical properties \cite{Rou53,SSS80,DoE94}. Here, we consider
a very crude approximation where a planar polymer is modeled as
a collection of $n$ rigid rods, with equal fixed lengths $l_{0}$ and we denote their
extremities by $\left(X_{0},X_{1},X_{2},\ldots,X_{n} \right)$ (Figure
\ref{figure1}) in a framework with origin ${\bf 0}$. We shall fix
one of the polymer ends $X_{0}=(L,0)$ (where $L>0$) on the $x$-axis.
The dynamics of the $i$-th rod is characterized by its angle
$\theta_{i}(t)$ with respect to the $x$-axis. The overall polymer
dynamics is thus characterized by the angles
$\left(\theta_{1}(t),\theta_{2}(t),\ldots,\theta_{n}(t), t\geq 0
\right)$. Due to the thermal collisions in the medium, each
angle follows a Brownian motion. Thus, with $\stackrel{(law)}{=}$ denoting equality in law,
\beqq
\left( \theta_{i}(t),i\leq n \right) \stackrel{(law)}{=} \sqrt{2D} \left( B_{i}(t),i\leq n \right) \ \mathrm{that} \ \mathrm{is:} \ \
d\theta_{i}(t) = \sqrt{2D} \: dB_{i}(t), \ \ i\leq n,
\eeqq
where $D$ is the rotational diffusion constant and $\left(
B_{1}(t),\ldots,B_{n}(t), t\geq 0 \right)$ is an $n$-dimensional
Brownian motion (BM). The position of each rod can now be obtained
as:
\beq\label{2Xi}
  X_{1}(t) = L + l_{0} e^{i\theta_{1}(t)}, \ X_{2}(t) = X_{1}(t) + l_{0} e^{i\theta_{2}(t)}, \ \ldots, \
  X_{n}(t) = X_{n-1}(t) + l_{0} e^{i\theta_{n}(t)}. \nonumber \\
\eeq
In particular, the moving end is given by:
\beq\label{2Xn}
  X_{n}(t) = L + l_{0} \sum^{n}_{k=1}   e^{i\theta_{k}(t)} = L + l_{0} \sum^{n}_{k=1}   e^{i\sqrt{2D} B_{k}(t)}, \nonumber \\
\eeq
which can be written as:
\beq\label{2XnRphi}
  X_{n}(t) = R_{n}(t) e^{i \varphi_{n}(t)}.
\eeq
Thus $\varphi_{n}(t)$ accounts for the rotation of the polymer with
respect to the origin ${\bf 0}$ and $R_{n}$ is the distance to the origin.

In order to compute the MRT, we shall study a sum of exponentials of
Brownian motions, a topic which often leads to surprising computations
\cite{Yor01}.
First, we scale the space and time variables as follows:
\beq\label{rescaling}
    \tilde{l}=\frac{L}{l_{0}}  \ \ \mathrm{and} \ \ \tilde{t}=\frac{t}{2D}.
\eeq
Equation (\ref{2Xn}) becomes:
\beq\label{2Xnrescaled}
  X_{n}(t) = \tilde{l} + \sum^{n}_{k=1}   e^{i\tilde{B}_{k}(t)},
\eeq
where $\left(\tilde{B}_{1}(t),\ldots,\tilde{B}_{n}(t), t\geq 0 \right)$ is an $n$-dimensional
Brownian motion (BM), and for $k=1,\ldots,n$, using the scaling property of Brownian motion, we have:
\beqq
    \tilde{B}_{k}(t)\equiv \frac{1}{\sqrt{2D}} B_{k}(t) \stackrel{(law)}{=} B_{k}\left(\frac{t}{2D}\right)=B_{k}(\tilde{t}).
\eeqq

Before describing our approach, we first discuss
the mean initial configuration of the polymer. It is given by:
\beq\label{initialbis2}
c_{n}= E \left( \sum^{n}_{k=1} e^{i\theta_{k}(0)} \right),
\eeq
where the initial angles $\theta_{k}(0)$ are such that the polymer
has not already made a loop. After scaling, the mean initial
configuration becomes:
\beq\label{initial}
c_{n}= E \left( \sum^{n}_{k=1} e^{i \: \tilde{\theta}_{k}(0)} \right),
\eeq
with
\beq\label{rescalinginitialangle}
    \tilde{\theta}_{k}(0)=\tilde{B}_{k}(0).
\eeq
From now on, we use $\theta_{k}$ instead of
$\tilde{\theta}_{k}$, $B$ instead of $\tilde{B}$ and $t$ instead of $\tilde{t}$. Any segment in the interior of the
polymer can hit the angle $2\pi$ around the origin, but we will not
consider this as a winding event, although we could and in that
case, the MRT would be different.
Rather, we shall only consider that, given an initial configuration
$c_{n}$, the MRT is defined as:
\beq
MRT \equiv E\{\tau_{n}|c_{n} \} \equiv  E\left[\tau_{n}\right],
\eeq
where
\beq\label{deftaun}
\tau_{n} \equiv \inf\{ t>0, |\varphi_{n}(t)|=2\pi \}.
\eeq
Thus, an initial configuration is not winding when
\beq
\left| \varphi_{n}(0) \right| <2\pi.
\eeq
Then, we can define the winding event using a one dimensional
variable only. In general, winding is a rare event and we expect
that the MRT will depend crucially on the length of the polymer which
will be quite long. Interestingly, the rotation is accomplished
when the angle $\varphi_{n}(t)$ reaches $2 \pi$ or $-2 \pi$, but
the distance of the free end point to the origin is not fixed,
leading to a one dimensional free parameter space. This undefined
position is in favor of a winding time that is not too large when
compared to any narrow escape problem where a Brownian particle has
to find a small target in a confined domain
\cite{WaK93,HoS04,SSHE06,SSH07,BKH07,PWPSK09}.

In this study, we consider not only that the initial
condition satisfies $\left| \varphi_{n}(0) \right| <2\pi$, but we
impose that $\varphi_{n}(0)$ is located far enough from
$2\pi$, to avoid studying any boundary layer effect, which would
lead to a different MRT law. Indeed, starting inside the boundary
layer for a narrow escape type problem leads to specific escape laws
\cite{SSHE06}. Given a small $\varepsilon>0$ we shall consider the space of configurations $\Omega_{\varepsilon}$
such that  $\left|\varphi_{n}(0) \right| <2\pi-\varepsilon$. We shall mainly focus
on the stretched polymer,
\beq
\left(\theta_{1}(0),\theta_{2}(0),\ldots,\theta_{n}(0)
\right)=\left(0,0,\ldots,0\right),
\eeq
and thus $c_{n}=n >0$ (in this case, $\varphi_{n}(0)=0$).
Finally, it is quite obvious that winding occurs only when the
condition:
\beq
n l_{0}>L
\eeq
is satisfied, which we assume all along.

The outline of our method is: first we show that the sum $X_{n}(t)$
converges  (eq. (\ref{2Xn})) and we obtain a Central Limit Theorem.
Using It\^{o} calculus, we study the sequence:
\beq
\frac{1}{\sqrt{n}} \overline{X}_{n}(t) = \frac{1}{\sqrt{n}} \left[ X_{n}(t) - E\left( X_{n}(t) \right) \right].
\eeq
and prove that $\frac{1}{\sqrt{n}} \overline{X}_{n}(t)$, for $n$ large, converges to a stochastic
process which is a generalization of an Ornstein-Uhlenbeck process
(GOUP), containing a time dependent deterministic drift
$c_{n}e^{-t}$. This GOUP is driven by a martingale
$\left(M_{t}^{(n)}, t\geq 0 \right)$ that we further characterize.
Interestingly, the two cartesian coordinates of
$\left(M_{t}^{(n)}, t\geq 0 \right)$ converge to two independent
Brownian motions with two different time scale functions. To obtain
an asymptotic formula for the MRT, we show that in the long time
asymptotics, where winding occurs, the GOUP can be approximated by a
standard Ornstein-Uhlenbeck process (OUP). Using some properties of
the GOUP \cite{Vak10}, we finally derive the MRT for the polymer
which is the mean time that $|\varphi_{n}(t)|=2\pi$.

\section{Properties of the free polymer end $X_{n}(t)$ using a Central Limit Theorem}
In this section we study some properties of the free polymer end
$X_{n}(t)$. In particular, using a Central Limit Theorem,
we show that the limit process satisfies a stochastic equation of a
new type. To study the random part of $X_{n}(t)$, we shall remove
from it its first moment and we shall now consider the asymptotic
behavior of the drift-less sequence:
\beq \label{2Xnbar1}
    \frac{1}{\sqrt{n}} \overline{X}_{n}(t) =
    \frac{1}{\sqrt{n}} \left[ X_{n}(t) - E\left( X_{n}(t) \right) \right].
\eeq
We start by computing the first moment $E\left( X_{n}(t) \right)$.
Because (see (\ref{2Xn})), $\left(\theta_{i}(t),t\geq0\right)$
are assumed to be $n$ independent identically distributed (iid) Brownian motions
with variance $2D$, after rescaling, we obtain:
\beq \label{meanX_{n}(t)}
    E\left( X_{n}(t) \right) &=& E \left[\tilde{l} + \sum^{n}_{k=1} e^{i B_{k}(t)}\right] \nonumber \\
    &=& \tilde{l} + \left(\sum^{n}_{k=1}  E \left[e^{i\left(B_{k}(t)-B_{k}(0)\right)}\right] \: E \left[e^{i\left(B_{k}(0)\right)}\right]\right) \nonumber \\
    &=& \tilde{l} + c_{n} e^{-\frac{t}{2}}.
\eeq
where  $c_{n}$ is defined by (\ref{initial}) and we have used that:
\beq
 E \left[e^{i\left(B_{k}(t)-B_{k}(0)\right)}\right]= e^{-\frac{t}{2}} \ .
\eeq
We study the sequence (\ref{2Xnbar1}) as follows:
\beq \label{2Xnbar}
    \frac{1}{\sqrt{n}} \overline{X}_{n}(t) &=& \frac{1}{\sqrt{n}} \left[ \sum^{n}_{k=1} e^{iB_{k}(t)} - E\left( \sum^{n}_{k=1} e^{i B_{k}(t)} \right) \right] \nonumber \\
    &=& \frac{1}{\sqrt{n}} \sum^{n}_{k=1} F_{k} (t),
\eeq
where $F_{k} (t) = e^{iB_{k}(t)} - E\left(e^{iB_{k}(t)} \right)$. Applying It\^{o}'s formula to
\beq\label{ZsumF}
Z_{t}^{(n)} =
\frac{1}{\sqrt{n}} \sum^{n}_{k=1} F_{k} (t) \ ,
\eeq
with $Z_{0}^{(n)}=0$, we obtain:
\beq
Z_{t}^{(n)} &=& \frac{i}{\sqrt{n}} \int^{t}_{0} \sum^{n}_{k=1}
e^{i B_{k}(s)} dB_{k}(s) - \frac{1}{2\sqrt{n}} \int^{t}_{0}
\sum^{n}_{k=1} \left( e^{i B_{k}(s)} - E\left(
e^{i B_{k}(s)} \right)
\right) ds \nonumber \\
\label{2Zn}
    &=& M_{t}^{(n)} - \frac{1}{2}\int^{t}_{0} Z_{s}^{(n)} \: ds,
\eeq
where:
\beq\label{2MSC}
M_{t}^{(n)} &=& \frac{i}{\sqrt{n}} \int^{t}_{0} \sum^{n}_{k=1} e^{i B_{k}(s)} \: dB_{k}(s) \nonumber \\
    &=& \frac{1}{\sqrt{n}} \int^{t}_{0} \sum^{n}_{k=1} \left( i\cos(B_{k}(s))- \sin( B_{k}(s)) \right) \: dB_{k}(s) \nonumber \\
    &=& -S_{t}^{(n)} + i C_{t}^{(n)}.
\eeq
We shall now study the asymptotic limit of the martingales
$M_{t}^{(n)}$ as $n\rightarrow\infty$ and summarize our result in the following theorem;
the convergence in law being considered there is associated with the topology of the
uniform convergence on compact sets of the functions in
$C(\mathbb{R}_{+},\mathbb{R}^{2})$.

\begin{theo}\label{CLT}
The sequence $(M_{t}^{(n)},t\geq 0)$ converges in law to a Brownian
motion in dimension 2, with 2 deterministic time changes. More precisely,
\beq \label{2CLTcvg}
    (S^{(n)}_{t},C^{(n)}_{t}, t\geq0) \overset{{(law)}}{\underset{n\rightarrow\infty}\longrightarrow} ( \sigma_{\left(\frac{1}{2}\int^{t}_{0} ds \: (1-e^{-2s})\right)} , \gamma_{\left(\frac{1}{2}\int^{t}_{0} ds \: (1+e^{-2s})\right)} , t\geq0),
\eeq
where $(\sigma_{u},\gamma_{u}, u\geq0)$ are two independent Brownian motions.
\end{theo}
\begin{rem}\label{remTCL}
$\left.a\right)$ The convergence in law (\ref{2CLTcvg}) is a new result and it
shows that the sum of the complex exponentials of i.i.d. Brownian
motions can be approximated by a two dimensional Brownian motion,
with a different time scale for each coordinate. \\
$\left.b\right)$ An extension of Theorem \ref{CLT} when $\left(\exp\left(iB(s)\right),s\geq0\right)$
a Brownian motion on the unit circle is replaced by a BM on the unit sphere in $\mathbb{R}^{n}$
is obtained in \cite{HVY11}.
\end{rem}

{\noindent \textbf{Proof of Theorem \ref{CLT}}} See Appendix \ref{apTCL}. \hfill \QED

\vspace{8pt}
\noindent We conclude that the sequence $(M_{t}^{(n)},t\geq 0) $
converges in law:
\beq\label{2Mcvg}
    (M_{t}^{(n)},t\geq 0) \overset{{(law)}}{\underset{n\rightarrow\infty}\longrightarrow} (\sigma_{\left(\frac{t}{2}-\frac{1-e^{-2t}}{4}\right)} + i \gamma_{\left(\frac{t}{2}+\frac{1-e^{-2t}}{4}\right)},t\geq 0).
\eeq
The process $(Z_{t}^{(n)},t\geq 0)$, which we defined in (\ref{2Zn}), is a generalization of the classical
Ornstein-Uhlenbeck process; it is driven by $(M_{t}^{(n)},t\geq 0)$
and, from (\ref{2Zn}), we obtain:
\beq\label{2OU}
    Z_{t}^{(n)} &=& e^{-\frac{t}{2}} \int^{t}_{0} e^{\frac{s}{2}} dM_{s}^{(n)}.
\eeq

\begin{corr}\label{2Corr}
$\left.\mathrm{a}\right)$ The sequence $(Z_{t}^{(n)},t\geq 0)$ converges in law:
\beq
Z_{t}^{(n)}\overset{{(law)}}{\underset{n\rightarrow\infty}\longrightarrow}Z_{t}^{(\infty)},
\eeq
with:
\beq\label{cvgZn}
    Z_{t}^{(\infty)} = e^{-\frac{t}{2}} \int^{t}_{0} \sqrt{\sinh(s)} \; d\delta_{s} + i \; e^{-\frac{t}{2}} \int^{t}_{0} \sqrt{\cosh(s)} \; d\tilde{\delta}_{s},
\eeq
where $\left(\delta_{t}, \tilde{\delta}_{t},t\geq0\right)$ are two independent 1-dimensional Brownian motions.
$\left.\mathrm{b}\right)$ $Z_{t}^{(\infty)}\overset{{(law)}}{\underset{t\rightarrow\infty}\longrightarrow}\frac{1}{\sqrt{2}}\left(N+i\tilde{N}\right)$, where $N$ and $\tilde{N}$ are two centered and reduced Gaussian variables.
\end{corr}

{\noindent \textbf{Proof of Corollary \ref{2Corr}}} \\
$\left.\mathrm{a}\right)$ From (\ref{2MSC}), (\ref{2CLTcvg}) and (\ref{2OU}) we deduce:
\beq
Z_{t}^{(n)} &=& e^{-\frac{t}{2}} \int^{t}_{0} e^{\frac{s}{2}} dM_{s}^{(n)}\label{2Zncvg} \\
    &\overset{{(law)}}{\underset{n\rightarrow\infty}\longrightarrow}& e^{-\frac{t}{2}} \int^{t}_{0} \underbrace{e^{\frac{s}{2}} \sqrt{\frac{1-e^{-2s}}{2}}}_{\sqrt{\sinh(s)}} \; d\delta_{s}
    + i e^{-\frac{t}{2}} \int^{t}_{0} \underbrace{e^{\frac{s}{2}} \sqrt{\frac{1+e^{-2s}}{2}}}_{\sqrt{\cosh(s)}} \; d\tilde{\delta}_{s}  \equiv Z_{t}^{(\infty)}. \nonumber \\
     \label{2Zncvgfinal}
\eeq
$\left.\mathrm{b}\right)$ From (\ref{2Zncvgfinal}), we change variables: $u=t-s$ and we obtain:
\beq
    Z_{t}^{(\infty)} &=& \int^{t}_{0} e^{-\frac{t-s}{2}} \sqrt{\frac{1-e^{-2s}}{2}} \; d\delta_{s}
    + i \int^{t}_{0} e^{-\frac{t-s}{2}} \sqrt{\frac{1+e^{-2s}}{2}} \; d\tilde{\delta}_{s} \nonumber \\
    &\overset{{u=t-s}}{\underset{(law)}=}& \int^{t}_{0} e^{-\frac{u}{2}} \sqrt{\frac{1}{2}-\frac{e^{-2(t-u)}}{2}} \; d\delta_{u}
    + i \int^{t}_{0} e^{-\frac{u}{2}} \sqrt{\frac{1}{2}+\frac{e^{-2(t-u)}}{2}} \; d\tilde{\delta}_{u} \nonumber \\
    && \label{2Zncvgfinalbis} \\
    &\overset{{L^{2}}}{\underset{t\rightarrow\infty}\longrightarrow}& \frac{1}{\sqrt{2}} \int^{\infty}_{0} e^{-\frac{u}{2}}d\delta_{u} + i \; \frac{1}{\sqrt{2}} \int^{\infty}_{0} e^{-\frac{u}{2}} d\tilde{\delta}_{u} \ , \label{2Zncvgfinalbisbis}
\eeq
where the two variables on the RHS of (\ref{2Zncvgfinalbisbis}) are centered Gaussian with variance $1/2$ and the convergence in $L^{2}$ for $t\rightarrow\infty$ may be proved by using the dominated convergence theorem. \hfill \QED
\begin{rem}\label{2ItoOchi}
In \cite{It\^{o}83} and \cite{Ochi85} there is a more complete asymptotic study for $Z_{t}^{(n)}(\varphi)$,
as $n\rightarrow\infty$, for a large class of functions $\varphi$, but we do not pursue this line of research here.
\end{rem}
On one hand, from the identities (\ref{2Xnbar1}), (\ref{meanX_{n}(t)}),
(\ref{2Xnbar}) and (\ref{ZsumF}), we obtain the following expansion:
\beq\label{2Xnfinal}
    X_{n}(t) &=& \overline{X}_{n}(t) + E\left[ X_{n}(t) \right] \nonumber \\
             &=& \sqrt{n} \: Z_{t}^{(n)} + c_{n} e^{-\frac{t}{2}} + \tilde{l} \nonumber \\
             \Rightarrow \tilde{X}_{n}(t)&\equiv& X_{n}(t)-\tilde{l}=\sqrt{n} \: Z_{t}^{(n)} + c_{n} e^{-\frac{t}{2}} ,
\eeq
with $n$ the number of rods/beads, $Z_{t}^{(n)}$ a GOUP driven by $M_{t}^{(n)}$ which is given by
(\ref{2Mcvg}), $c_{n}=E \left( \sum^{n}_{k=1} e^{i\theta_{k}(0)} \right)$ a constant depending on the mean initial
configuration and $\tilde{l}=\frac{L}{l_{0}}$ the rescaled distance of the fixed end from
the origin ${\bf 0}$ ($l_{0}$ is the fixed length of the rods). \\
On the other hand, with Corollary \ref{2Corr}(a), we deduce that:
\beq \label{2Xnbar2}
    Z^{(n)}_{t} &\equiv& \frac{1}{\sqrt{n}} \overline{X}_{n}(t) \equiv \frac{1}{\sqrt{n}} \left[ X_{n}(t) - E\left( X_{n}(t) \right) \right] \nonumber \\
    &=& \frac{1}{\sqrt{n}} \left[ X_{n}(t) - c_{n} e^{-\frac{t}{2}} - \tilde{l}\right] \nonumber \\
    &\overset{{(law)}}{\underset{n\rightarrow\infty}\longrightarrow}& e^{-\frac{t}{2}} \int^{t}_{0} \sqrt{\sinh(s)} \; d\delta_{s} + i \; e^{-\frac{t}{2}} \int^{t}_{0} \sqrt{\cosh(s)} \; d\tilde{\delta}_{s}\equiv Z^{(\infty)}_{t}.
\eeq

\section{Asymptotic expression for the MRT}\label{2secMRT}
We now study more precisely the different time scales of the
two Brownian motions in (\ref{2Mcvg}) and (\ref{2Xnbar2}). We estimate $Z_{t}^{(\infty)}$ for $t$ large,
the regime for which  the rotation will be accomplished. \\
We introduce here the following notation: $\stackrel{L^{2}}{\approx}$ denotes
closeness in the $L^{2}-$norm: for two stochastic processes $(W^{(1)}_{t},t\geq0)$ and
$(W^{(2)}_{t},t\geq0)$, the notation $W^{(1)}_{t}\stackrel{L^{2}}{\approx}W^{(2)}_{t}$ means that
${\underset{t\rightarrow \infty}\lim} E\left[\left|W^{(1)}_{t}-W^{(2)}_{t}\right|^{2}\right]=0$. \\
We shall show that, with $\left(\mathbb{B}_{t} = \delta_{t}+i
\tilde{\delta}_{t},t\geq0\right)$ a 2-dimensional Brownian
motion starting from $1$:
\beq\label{2I}
  Z_{t}^{(\infty)} &=& e^{-\frac{t}{2}} \int^{t}_{0} \sqrt{\sinh(s)} \; d\delta_{s} + i \; e^{-\frac{t}{2}} \int^{t}_{0} \sqrt{\cosh(s)} \; d\tilde{\delta}_{s} \nonumber \\
  &\stackrel{L^{2}}{\approx}& e^{-\frac{t}{2}} \int^{t}_{0} \frac{e^{\frac{s}{2}}}{\sqrt{2}} d\mathbb{B}_{s}.
\eeq
For this, it suffices to use the expression (\ref{2Zncvgfinalbisbis})
and the following Proposition, which reinforces the $\stackrel{L^{2}}{\approx}$
result in (\ref{2I}).
\begin{prop}\label{propositionZ}
As $t\rightarrow\infty$, the Gaussian martingales
\beqq
    \left(\int^{t}_{0} \sqrt{\sinh(s)} \; d\delta_{s}-\int^{t}_{0} \frac{e^{s/2}}{\sqrt{2}} \; d\delta_{s}, t\geq0\right),
\eeqq
and
\beqq
    \left(\int^{t}_{0} \sqrt{\cosh(s)} \; d\tilde{\delta}_{s}-\int^{t}_{0} \frac{e^{s/2}}{\sqrt{2}} \; d\tilde{\delta}_{s}, t\geq0\right)
\eeqq
converge a.s. and in $L^{2}$. The limit variables are Gaussian with variances $\frac{\pi-3}{2}$ and $-1+2\sqrt{2}-2a_{s}(1)\approx0,033$, where $a_{s}(x)\equiv\arg\sinh (x) \equiv \log (x+\sqrt{1+x^{2}}), \ x \in \mathbb{R}$, respectively.
\end{prop}
Thus, by multiplying both processes by $e^{-\frac{t}{2}}$, we obtain (\ref{2I}).\\
\vspace{5pt}

{\noindent \textbf{Proof of Proposition \ref{propositionZ}}}
The Gaussian martingale $\int^{t}_{0} \left(\sqrt{\sinh(s)} - \frac{e^{s/2}}{\sqrt{2}}\right) \; d\delta_{s}$ has
increasing process
\beqq
    \int^{t}_{0} \left(\sqrt{\sinh(s)}-\frac{e^{s/2}}{\sqrt{2}}\right)^{2} \; ds = \int^{t}_{0} \frac{e^{s}}{2}\left(\sqrt{1-e^{-2s}}-1 \right)^{2} \; ds \ ,
\eeqq
which converges as $t\rightarrow\infty$. Hence, the limit variable $\int^{\infty}_{0} \left(\sqrt{\sinh(s)} - \frac{e^{s/2}}{\sqrt{2}}\right) \; d\delta_{s}$ is Gaussian, and its variance is given by (we change variables: $u=e^{-2s}$ and $B(a,b)$ denotes the Beta function with arguments $a$ and $b$\footnote[7]{We recall that if $\left(\Gamma(x),x\geq0\right)$ denotes the Gamma function, then $B(a,b)=\frac{\Gamma(a)\Gamma(b)}{\Gamma(a+b)}$.}):
\beqq
    && \int^{\infty}_{0} \frac{ds \; e^{s}}{2} \left(\sqrt{1-e^{-2s}} - 1\right)^{2} = \frac{1}{4} \int^{1}_{0} du \; u^{-3/2} \left(\sqrt{1-u} - 1\right)^{2} \\
    &=& \frac{1}{4} \left[ \int^{1}_{0} du \; u^{-3/2} \left((1-u)-2\sqrt{1-u}+1\right)\right] \\
    &=& \frac{1}{4} \left\{B\left(-\frac{1}{2},2\right)-2B\left(-\frac{1}{2},\frac{3}{2}\right)-2\right\} = \frac{\pi-3}{2} \ .
\eeqq
To be rigorous, the integral $\int^{1}_{0} du \; u^{-\alpha} \left(\sqrt{1-u} - 1\right)^{2}$, which is well defined for $0<\alpha<1$,
can be extended analytically for any complex $\alpha$ with $\mathrm{Re}(\alpha)<3$. \\
For the convergence of the second process, it suffices to replace $\sinh(s)$ by $\cosh(s)\equiv\frac{e^{s}}{2}(1+e^{-2s})$.
The limit variable $\int^{\infty}_{0} \left(\sqrt{\cosh(s)} - \frac{e^{s/2}}{\sqrt{2}}\right) \; d\tilde{\delta}_{s}$ is also Gaussian,
and, repeating the previous calculation, we easily compute its variance. \hfill \QED

\subsubsection*{Asymptotic expression for the MRT}\label{2asymptMRT}
For $t$ large, we derive an asymptotic value for the MRT.
First, from (\ref{2OU}) and (\ref{2Xnfinal}):
\beq\label{2Xnfinal3}
    \tilde{X}_{n}(t) &=& \sqrt{n}Z^{(n)}_{t}+c_{n}e^{-\frac{t}{2}} \nonumber \\
    &=& \sqrt{n}e^{-\frac{t}{2}} \left( \tilde{c}_{n} + \int^{t}_{0} e^{\frac{s}{2}} dM^{(n)}_{s} \right),
\eeq
where, the sequence $\tilde{c}_{n}$ is:
\beq\label{2cn'}
\tilde{c}_{n}\equiv \frac{c_{n}}{\sqrt{n}},
\eeq
$n$ is the number of rods/beads, and $c_{n} \equiv
E \left( \sum^{n}_{k=1} e^{i \theta_{k}(0)} \right)$ is a constant depending on the
mean initial configuration. Thus, from (\ref{2Xnfinal}), (\ref{2I}), (\ref{2Xnfinal3}) and also using the scaling property of Brownian motion:
\beq\label{2Xnfinal3bis}
    \tilde{X}_{n}(t) \overset{{(law)}}{\underset{n:large}\approx} \sqrt{n}Y_{t}^{(n)},
\eeq
where
\beq\label{2Xnfinal4}
    Y_{t}^{(n)}\equiv e^{-\frac{t}{2}} \left( \tilde{c}_{n} + \int^{t}_{0} e^{\frac{s}{2}} d\mathbb{B}_{s/2} \right).
\eeq

\subsubsection*{Changing time and expression of the MRT}\label{2asymptMRTpp}
To express the MRT, we now apply deterministic time changes. To make our writing simple, we denote $Y_{t}$ for $Y_{t}^{(n)}$,
and changing variables $u=\frac{s}{2}$ in (\ref{2Xnfinal4}), we obtain:
\beq\label{2OUz0DeB2}
    Y_{2t} = e^{-t} \left( \tilde{c}_{n} +  \int^{t}_{0} e^{u} d\mathbb{B}_{u} \right).
\eeq
Now, there is another BM $(\tilde{\mathbb{B}}_{t},t\geq0)$, starting from $\tilde{c}_{n}$, such that:
\beq\label{2OUz0DeBDS}
    Y_{2t} = e^{-t} \left( \tilde{\mathbb{B}}_{\alpha_{t}} \right),
\eeq
where:
\beqq
     \alpha_{t}= \int^{t}_{0} e^{2s} ds =  \frac{e^{2t}-1}{2},
\eeqq
hence:
\beq
     \alpha^{-1}(t)= \frac{1}{2}  \ln \left( 1+2t\right).\label{eeeeq}
\eeq
Applying It\^{o}'s formula to (\ref{2OUz0DeBDS}), we
obtain:
\beqq
dY_{2s} = -e^{-s} \tilde{\mathbb{B}}_{\alpha_{s}} ds + e^{-s} d\left(\tilde{\mathbb{B}}_{\alpha_{s}}\right).
\eeqq
We divide by $Y_{2s}$ and we obtain:
\beqq
\frac{dY_{2s}}{Y_{2s}} = - ds +\frac{d\tilde{\mathbb{B}}_{\alpha_{s}}}{\tilde{\mathbb{B}}_{\alpha_{s}}}.
\eeqq
Thus:
\beqq
\mathrm{Im} \left(\frac{dY_{2s}}{Y_{2s}}\right) =
\mathrm{Im}
\left(\frac{d\tilde{\mathbb{B}}_{\alpha_{s}}}{\tilde{\mathbb{B}}_{\alpha_{s}}}\right),
\eeqq
which means that, if we denote:
\beq
\theta^{Z}_{t}\equiv
\mathrm{Im}(\int^{t}_{0}\frac{dZ_{s}}{Z_{s}}), t\geq0,
\eeq
the continuous winding process associated to a generic stochastic
process $Z$, then:
\beqq
\theta^{Y}_{2t} =
\theta^{\tilde{\mathbb{B}}}_{\alpha_{t}}.
\eeqq
Thus, the first hitting times of the symmetric conic
boundary of angle $c$
\beq\label{2Tthetac}
T^{|\theta^{Y}|}_{c}\equiv
\inf \left\{t\geq 0 : \left|\theta^{Y}_{t}\right|=c \right\},
\eeq
and
\beq\label{2Tthetacc}
T^{|\theta^{\tilde{\mathbb{B}}}|}_{c}\equiv\inf \left\{t\geq 0 : \left|\theta^{\tilde{\mathbb{B}}}_{t}\right|=c \right\},
\eeq
for an Ornstein-Uhlenbeck process $Y$ with parameter $\lambda=1$, and for a Brownian motion $\tilde{\mathbb{B}}$
respectively, with relation (\ref{eeeeq}), satisfy:
\beq\label{2bisTchat}
    2T^{|\theta^{Y}|}_{c}=\frac{1}{2}\ln
    \left(1+2T^{|\theta^{\tilde{\mathbb{B}}}|}_{c}\right).
\eeq
Finally,
\beq\label{ETchatz0Dasympbis}
    E\left[2T^{|\theta^{Y}|}_{c}\right] &=& \frac{1}{2} E\left[\ln \left(1+2 T^{|\theta^{\tilde{\mathbb{B}}}|}_{c}\right) \right] \\
                        &=& \frac{\ln 2}{2} + \frac{1}{2} E\left[\ln\left(T^{|\theta^{\tilde{\mathbb{B}}}|}_{c}+\frac{1}{2}\right)\right] ,
\eeq
and equivalently:
\beq\label{ETchatz0Dasymp2}
    E\left[T^{|\theta^{Y}|}_{c}\right] &=& \frac{\ln 2}{4} + \frac{1}{4} E\left[\ln\left(T^{|\theta^{\tilde{\mathbb{B}}}|}_{c}+\frac{1}{2}\right)\right].
\eeq
Thus, by taking $c=2\pi$, for $n$ large, with $\tilde{\varphi}_{n}(t)$ denoting the total angle of $\tilde{X}_{n}(t)$,
the mean time $E\left[\tilde{\tau}_{n}\right]$, where $\tilde{\tau}_{n}\equiv \inf\{ t>0, |\tilde{\varphi}_{n}(t)|=2\pi \}$,
that $\tilde{X}_{n}(t)$ rotates around ${\bf 0}$, is:
\beq\label{2taunapprox}
E\left[\tilde{\tau}_{n}\right] &\approx&  \frac{\sqrt{n} }{4} \left( \ln 2 + E\left[\ln\left(T^{|\theta^{\tilde{\mathbb{B}}}|}_{2\pi}+\frac{1}{2}\right)\right]\right).
\eeq
By using the series expansion of $\log(1+x)$, we obtain informally:
\beq
    E\left[\ln\left(T^{|\theta^{\tilde{\mathbb{B}}}|}_{c}+\frac{1}{2}\right)\right] - E\left[\ln\left(T^{|\theta^{\tilde{\mathbb{B}}}|}_{c}\right)\right] &=& E\left[\ln\left(1+\frac{1}{2T^{|\theta^{\tilde{\mathbb{B}}}|}_{c}}\right)\right] \label{2logseries2} \\
    &=& \frac{1}{2} E\left[\frac{1}{T^{|\theta^{\tilde{\mathbb{B}}}|}_{c}}\right]-\frac{1}{8}E\left[\left(\frac{1}{T^{|\theta^{\tilde{\mathbb{B}}}|}_{c}}\right)^{2}\right]+\ldots, \nonumber \\ \label{2logseries}
\eeq
where this expansion should be understood as follows: we obtain an alternate series of negative moments of $T^{|\theta^{\tilde{\mathbb{B}}}|}_{c}$,
which converges \cite{VaY11}. For our purpose, we will truncate the series.  We shall estimate and use only the first moment of $1/T^{|\theta^{\tilde{\mathbb{B}}}|}_{c}$. By numerical calculations, the other moments can be neglected.

\begin{figure}
    \includegraphics[width=1.00\textwidth]{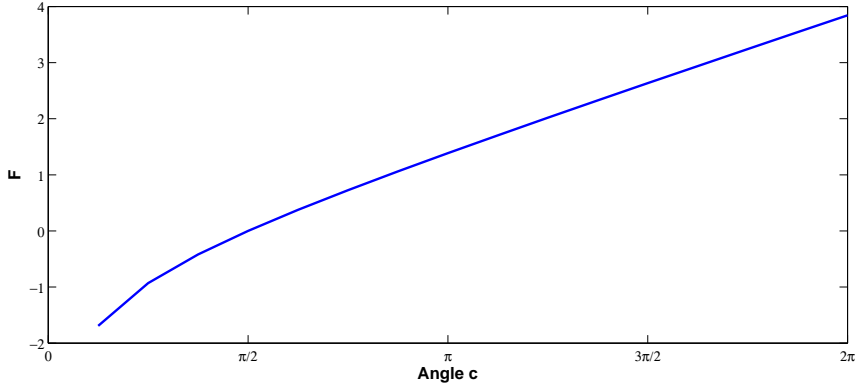} \\
    \caption{{\bf $F$ as a function of the angle $c$.} }\label{figFc}
\end{figure}

From the skew product representation \cite{KaSh88,ReY99} there is another planar Brownian motion $(\beta_{u}+i \gamma_{u},u\geq0)$, starting from $\log
\tilde{c}_{n}+i0$, such that:
\beq\label{2skew-product}
\log\left|\tilde{\mathbb{B}}_{t}\right|+i\theta_{t}\equiv\int^{t}_{0}\frac{d\tilde{\mathbb{B}}_{s}}{\tilde{\mathbb{B}}_{s}}=\left(
\beta_{u}+i\gamma_{u}\right)
\Bigm|_{u=H_{t}\equiv\int^{t}_{0}\frac{ds}{\left|\tilde{\mathbb{B}}_{s}\right|^{2}}},
\eeq
and equivalently:
\beq\label{2skew-product2}
\log\left|\tilde{\mathbb{B}}_{t}\right|=\beta_{H_{t}}; \ \ \theta_{t}=\gamma_{H_{t}}.
\eeq
Thus, with $T^{|\gamma|}_{c}\equiv\inf \left\{t\geq 0 :\left|\gamma_{t}\right|=c \right\}$ and
$T^{|\theta^{\tilde{\mathbb{B}}}|}_{c}\equiv\inf \left\{t\geq 0 :\left|\theta^{\tilde{\mathbb{B}}}_{t}\right|=c \right\}$,
because $\theta_{T^{|\theta^{\tilde{\mathbb{B}}}|}_{c}}=\gamma_{H_{T^{|\theta^{\tilde{\mathbb{B}}}|}_{c}}}$:
\beqq
T^{|\gamma|}_{c}=H_{T^{|\theta^{\tilde{\mathbb{B}}}|}_{c}},
\eeqq
hence $T^{|\theta^{\tilde{\mathbb{B}}}|}_{c}=H^{-1}_{u}\Bigm|_{u=T^{|\gamma|}_{c}}$, where:
\beq\label{2Hinverse}
H^{-1}_{u}\equiv \inf\{ t:H_{t}>u \} =
\int^{u}_{0}ds \exp(2\beta_{s})\equiv A_{u},
\eeq
and for $u=T^{|\gamma|}_{c}$, we obtain:
\beq\label{2negmom}
    T^{|\theta^{\tilde{\mathbb{B}}}|}_{c}=A_{T^{|\gamma|}_{c}}.
\eeq
It is of interest for our purpose to look at the first negative moment of $A_{t}$ which has the following integral representation for any $t>0$
(\cite{Duf00}, \cite{D-MMY00}, p.49, Prop. 7, formula (15)):
\beq\label{2negmom2}
    E\left[\frac{1}{A_{t}}\right]=\int^{\infty}_{0} y e^{-y^{2}t/2} \coth\left(\frac{\pi}{2}y\right) \; dy .
\eeq
With $\beta_{s}=\log(\tilde{c}_{n})+\beta^{(0)}_{s}$, where $(\beta^{(0)}_{s},s\geq0)$ is
a one-dimensional Brownian motion starting from 0, we obtain:
\beq \label{Tcthetacn}
T^{|\theta^{\tilde{\mathbb{B}}}|}_{c}&=& H^{-1}(T^{|\gamma|}_{c}) \equiv \int^{T^{|\gamma|}_{c}}_{0} ds \; \exp\left(2 \beta_{s}\right) \nonumber \\
&=& (\tilde{c}_{n})^{2} \left(\int^{T^{|\gamma|}_{c}}_{0} ds \; \exp\left(2 \beta^{(0)}_{s}\right)\right) \nonumber \\
&\equiv& (\tilde{c}_{n})^{2}T^{|\theta^{\mathbb{B}^{(1)}}|}_{c} \ \ ,
\eeq
where $T^{|\theta^{\mathbb{B}^{(1)}}|}_{c}\equiv \inf \left\{t\geq 0 :
\left|\theta^{\mathbb{B}^{(1)}}_{t}\right|=c \right\}$ is the first hitting time of the
symmetric conic boundary of angle $c$ of a Brownian motion $\mathbb{B}^{(1)}$ starting from $1+i0$.
Hence, from (\ref{2negmom}-\ref{2negmom2}) and by using the Laplace transform for the hitting time
$T^{|\gamma|}_{c}$ \cite{ReY99} (Chapter II, Prop. 3.7) or \cite{PiY03}(p.298):
\beq
E\left[e^{-\frac{y^{2}}{2}T^{|\gamma|}_{c}}\right]=\frac{1}{\cosh(y c)},
\eeq
we get:
\beq\label{21/T}
    E\left[\frac{1}{T^{|\theta^{\mathbb{B}^{(1)}}|}_{c}}\right] &=& \int^{\infty}_{0} y E\left[e^{-\frac{y^{2}}{2}T^{|\gamma|}_{c}}\right] \coth\left(\frac{\pi}{2}y\right) \; dy \nonumber \\
    &=& \int^{\infty}_{0} \frac{y}{\cosh(y c)} \coth\left(\frac{\pi}{2}y\right) \; dy\equiv G(c) .
\eeq
For $c=2\pi$, we obtain the numerical result:
\beq\label{num}
 G(2\pi)\equiv E\left[\frac{1}{T^{|\theta^{\mathbb{B}^{(1)}}|}_{2\pi}}\right]\approx 0.167.
\eeq
Thus, with (\ref{Tcthetacn}),
\beq \label{ETchatz0Dbisb}
    E\left[\ln\left(T^{|\theta^{\tilde{\mathbb{B}}}|}_{c}\right)\right] &=& 2\ln(\tilde{c}_{n}) + E\left[\ln\left(T^{|\theta^{\mathbb{B}^{(1)}}|}_{c}\right)\right],
\eeq
and from (\ref{2logseries}), we have:
\beq\label{2log1/2approx}
E\left[\ln\left(T^{|\theta^{\tilde{\mathbb{B}}}|}_{2 \pi}+\frac{1}{2}\right)\right] \approx 2\ln(\tilde{c}_{n}) + E\left[\ln\left(T^{|\theta^{\mathbb{B}^{(1)}}|}_{2\pi}\right)\right] + \frac{1}{2\tilde{c}^{2}_{n}}E\left[\frac{1}{T^{|\theta^{\mathbb{B}^{(1)}}|}_{2\pi}}\right]. \nonumber \\
\eeq
For the first moment of
$\ln\left(T^{|\theta^{\tilde{\mathbb{B}}}|}_{c}\right)$, for an angle $c$, using Bougerol's identity \cite{Bou83,ReY99}, there is
the integral representation \cite{CMY98,Vak10}:
\beq \label{ETchatz0Dbisb1}
E\left[\ln\left(T^{|\theta^{\mathbb{B}^{(1)}}|}_{c}\right)\right]=2F(c) + \ln\left(2\right) + c_{E},
\eeq
where\footnote[8]{We note that there is a simple relation between $F$ and $G$: $\frac{\pi^{2}}{4c} \ F'\left(\frac{\pi^{2}}{4c}\right)=c \ G(c)$.
For a more complete discussion, see \cite{Vak10}.}:
\beq\label{2F(c)}
F(c)=\int^{\infty}_{0} \frac{dz}{\cosh \left(\frac{\pi z}{2}\right)}
\ln\left(\sinh\left(cz\right)\right),
\eeq
and $c_{E}\approx 0.577$ denotes Euler's constant.\\
For $c=2\pi$, we have $F(2 \pi)\approx3.84$. In Figure \ref{figFc} we plot $F$ with respect to
the angle $c$.
\\
In summary, from (\ref{ETchatz0Dasymp2}), (\ref{num}), (\ref{2log1/2approx}) and
(\ref{ETchatz0Dbisb1}), we approximate $E\left[T^{|\theta^{Y}|}_{2\pi}\right]$:
\beq\label{2EThat}
E\left[T^{|\theta^{Y}|}_{2\pi}\right] \approx  \frac{1}{4} \left( 2\ln(\tilde{c}_{n}) +  Q +
\frac{1}{2\tilde{c}^{2}_{n}}E\left[\frac{1}{T^{|\theta^{\mathbb{B}^{(1)}}|}_{2\pi}}\right] \right),
\eeq
where:
\beq
Q=2F(2\pi) + 2\ln 2 + c_{E}
\eeq
is a constant with $F(2\pi)\approx3.84$, $c_{E}\approx 0.577$, and
$E\left[\frac{1}{T^{|\theta^{\mathbb{B}^{(1)}}|}_{2\pi}}\right]\approx
0.167$, thus:
\beq
Q\approx 9.54.
\eeq
Thus, by taking $c=2\pi$, for $n$ large, the mean time $E\left[\tilde{\tau}_{n}\right]$ that
$\tilde{X}_{n}(t)$ rotates around ${\bf 0}$, is:
\beq\label{2taunapproxbis}
E\left[\tilde{\tau}_{n}\right] &\approx&  \frac{\sqrt{n}}{4} \left( 2\ln(\tilde{c}_{n}) + \frac{0.08}{\tilde{c}^{2}_{n}}+ Q \right),
\eeq
and since:
\beq\label{2constantrandom}
\tilde{c}_{n}\equiv\frac{c_{n}}{\sqrt{n}},
\eeq
the MRT of $\tilde{X}_{n}(t)$ is given by the formula:
\beq\label{2ETfinalrandomIC}
 E\left[\tilde{\tau}_{n}\right] \approx \frac{\sqrt{n}}{4}  \left[ 2 \ln \left( \frac{c_{n}}{\sqrt{n}}\right)+  0.08\frac{n}{c^{2}_{n}}+ Q \right].
\eeq
For a long enough polymer, such that
$nl_{0}>>L\Rightarrow n>> \tilde{l}$, thus
$\tilde{l}=\frac{L}{l_{0}}$ is negligible with respect to $X_{n}(t)$
and from (\ref{2Xnfinal}) $X_{n}(t)\approx\tilde{X}_{n}(t)$.
For a mean initial configuration $c_{n}\equiv E \left( \sum^{n}_{k=1} e^{i
\theta_{k}(0)} \right)$, we obtain that the MRT of the polymer is given by the formula:
\beq\label{2ETfinalrandomIClong}
 E\left[\tau_{n}\right] \approx \frac{\sqrt{n}}{4} \left[ 2 \ln \left( \frac{E \left( \sum^{n}_{k=1} e^{i
\theta_{k}(0)}\right)}{\sqrt{n}} \right)+\frac{0.08 n}{\left(E \left( \sum^{n}_{k=1} e^{i\theta_{k}(0)}\right)\right)^{2}}+
Q \right].
\eeq
Finally,
using the unscaled variables (\ref{rescaling}), we obtain:
\beq\label{2ETfinaluncsaledrandomIClong}
 E\left[\tau_{n}\right] \approx \frac{\sqrt{n}}{8D} \left[ 2 \ln \left( \frac{E \left( \sum^{n}_{k=1} e^{\frac{i}{\sqrt{2D}}\theta_{k}(0)}\right)}{\sqrt{n}} \right)+ \frac{0.08 n}{\left(E \left( \sum^{n}_{k=1} e^{\frac{i}{\sqrt{2D}}\theta_{k}(0)}\right)\right)^{2}}+ Q \right]. \nonumber \\
\eeq
Expressions (\ref{2ETfinalrandomIClong}) and (\ref{2ETfinaluncsaledrandomIClong}) show that
the leading order term of the MRT depends on the initial
configuration, however this dependence is weak.

We consider now that the polymer is initially stretched $\left(\theta_{k}(0)=0, \ \forall k=(1,...,n)\right)$, hence $c_{n}=n$.
Thus, from (\ref{2ETfinaluncsaledrandomIClong}),
the MRT is approximately:
\beq\label{2ETfinaluncsaled}
 E\left[\tau_{n}\right] \approx \frac{\sqrt{n}}{8D}  \left[ \ln \left(n \right)+ 0.08\frac{1}{n}+ Q \right].
\eeq

\begin{figure}
\centering
    \includegraphics[width=0.60\textwidth]{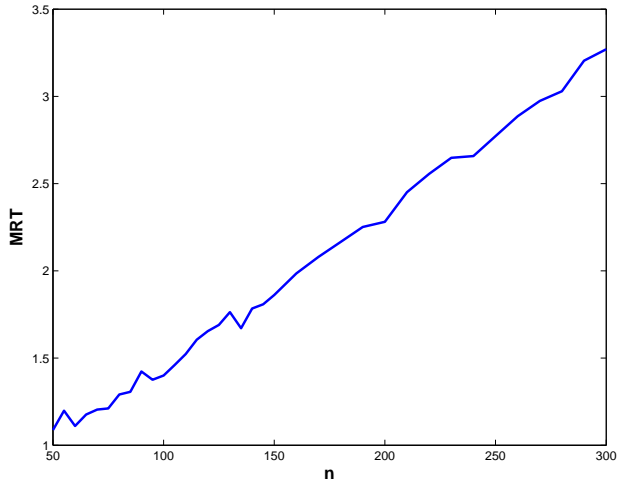} \\
    \caption{{\bf MRT of the free polymer end as a function of the number of beads $n$ (Brownian simulations)}.}\label{2E[T]}
\end{figure}

In order to check the range of validity of formula (\ref{2ETfinaluncsaled}), we ran some Brownian simulations. In Figure
\ref{2E[T]}, we simulated the MRT with a time step $dt=0.01$ for $n=50$ to $300$ rods in steps of 5,
and for each $n$, we took 300 samples and averaged over all of them.
The parameters we chose were $D=10$ for the diffusion coefficient,
$L=0.3$ for the distance from the origin ${\bf 0}$ and $l_{0}=0.25$ for the
length of each rod, and for the initial condition we chose a stretched
polymer located on the half line $\overrightarrow{{\bf 0}x}$, $\theta_{k}(0)=0$, $\forall k=(1,...,n)$
(hence $c_{n}=L+nl_{0}$), and then we computed the
MRT $E\left[\tau_{n}\right]$.
In Figure \ref{2comparisonimpr}, we plot both the results from
Brownian simulations and the formula (\ref{2ETfinaluncsaled}). We
considered values from $n=50$ to $300$ rods in steps of 10, $D=10$,
$L=0.3$, $l_{0}=0.25$, thus the condition $\sqrt{n} \: l_{0}>>L$
is satisfied. {For each numerical computation, we performed 300 runs
with a time step $dt=0.01$. By comparing the numerical
simulations (Fig. \ref{2comparisonimpr}), with the analytical formula for the  MRT, we see
an overshoot.
\begin{figure}
\centering
    \includegraphics[width=0.60\textwidth]{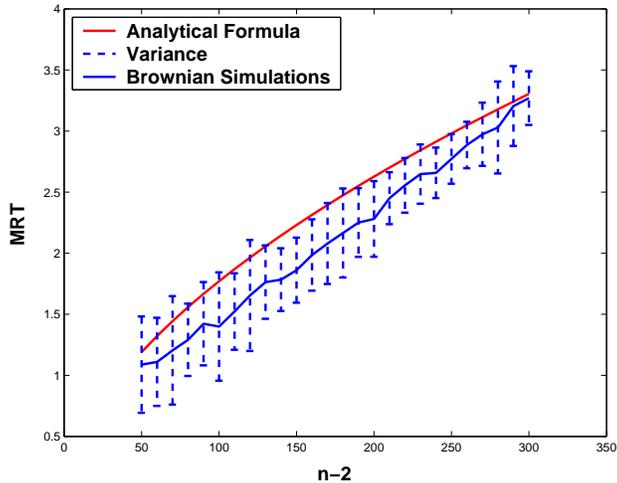} \\
    \caption{ {\bf Comparing the Analytical Formula and the Brownian simulations of the
MRT. } The MRT of the polymer end is depicted with respect to the
number of beads $n$. }\label{2comparisonimpr}
\end{figure}

The straight initial configuration is no restriction to the
generality of our study: we have $c_{n}\leq n$ and the upper
bound is achieved for a straight initial configuration. In Figure
\ref{straightrandom}, we simulate the MRT $E\left[\tau_{n}\right]$
for both a random and an initially straight configuration (Brownian
simulations with the same values for the parameters as above, after formula
(\ref{2ETfinaluncsaled})).
\begin{figure}
\centering
    \includegraphics[width=0.60\textwidth]{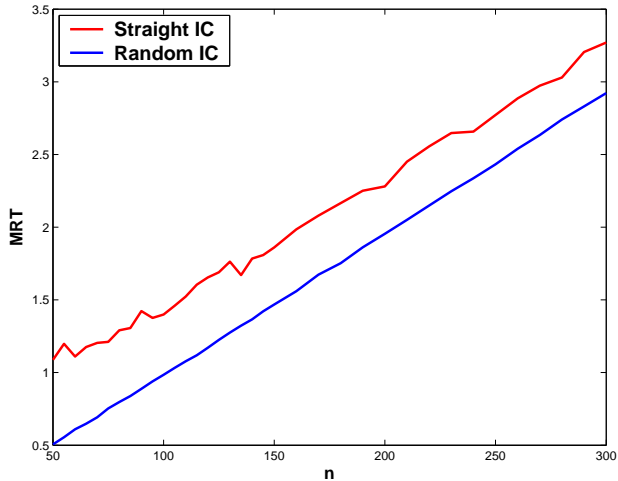} \\
    \caption{ {\bf MRT of the free polymer end as a function of the number of beads $n$
    for the straight and for a random initial configuration (Brownian simulations)}.}\label{straightrandom}
\end{figure}

\subsubsection*{Uniformly distributed initial angles}
When the initial angles
$\left(\theta_{1}(0),\theta_{2}(0),\ldots,\theta_{n}(0)\right)$ are
uniformly distributed over $\left[0,2\pi\right]$, by averaging over
all possible initial configurations, from (\ref{initial}) and
(\ref{rescalinginitialangle}), we obtain:
\beq
    c_{n}=E\left[\sum^{n}_{k=1} e^{i \theta_{k}(0)}\right]=E\left[\sum^{n}_{k=1} \left(\cos(\theta_{k}(0)) +i \; \sin(\theta_{k}(0))\right)\right]=0,
\eeq
hence, from (\ref{meanX_{n}(t)}),
\beq\label{2initialaverage}
    E\left[X_{n}(t)\right] \equiv \tilde{l}+c_{n} = \tilde{l}.
\eeq
We define:
\beq\label{XhattildeZ}
    \hat{X}_{n}(t)\equiv\frac{1}{\sqrt{n}} \left(X_{n}(t)-\tilde{l}\right)+1=1+Z^{(n)}_{t}.
\eeq
We know that, for $n$ large, with $\hat{\varphi}_{n}(t)$ denoting the total angle of $\hat{X}_{n}(t)$
the mean time $E\left[\hat{\tau}_{n}\right]$, where $\hat{\tau}_{n}\equiv \inf\{ t>0, |\hat{\varphi}_{n}(t)|=2\pi \}$,
that $\hat{X}_{n}(t)$ rotates around ${\bf 0}$, is:
\beqq
    E\left[\hat{\tau}_{n}\right]\approx  \frac{\sqrt{n}}{8D} \tilde{Q}.
\eeqq
Finally, for a long enough polymer, such that
$nl_{0}>>L\Rightarrow n>> \tilde{l}$ and $\sum^{n}_{k=1} e^{i \theta_{k}(t)}>>\sqrt{n}$, thus
$\tilde{l}=\frac{L}{l_{0}}$ and $\sqrt{n}$ are negligible with respect to $X_{n}(t)$
and from (\ref{XhattildeZ}) $X_{n}(t)\approx\hat{X}_{n}(t)$.
Using the unscaled variables (\ref{rescaling}), with $\tilde{Q}\approx 9.62$, the MRT satisfies:
\beq\label{2MRTrandominitial}
    E\left[\tau_{n}\right] \approx E\left[\hat{\tau}_{n}\right]\approx  \frac{\sqrt{n}}{8D} \tilde{Q}.
\eeq
\begin{rem}\label{2RealCase}
Formula (\ref{2ETfinaluncsaledrandomIClong}) or formulas (\ref{2ETfinaluncsaled}) and (\ref{2MRTrandominitial})
provide the asymptotic expansion for the MRT when  $\theta_{1}\in
\mathbb{R_{+}}$. In fact, $\theta_{1}$ is a reflected Brownian
motion in $[0,2\pi]$, thus a better
characterization is to estimate the MRT by using the probability
density function for $\theta_{1}$ in the one dimensional torus $[0,2\pi]$.
In the Appendix of \cite{Vakth11}, we derive this probability
density function and by repeating the previous calculations,
we show that, e.g. formula (\ref{2MRTrandominitial}), remains valid.
\end{rem}

\subsection{The Minimum Mean First Rotation Time}
The Minimum Mean Rotation Time (MMRT) is the first time that any of
the segments of the polymer loops around the origin,
\beq
MMRT \equiv \min_{\mathcal{E}_n} E\{\tau_{n}|c_{n} \}\equiv E\left[\tau_{min}\right],
\eeq
where $\mathcal{E}_n$ is the ensemble of rods which can travel up to
the origin. The MMRT is now a decreasing function of $n$. In Figure
\ref{figETvar}, we present some simulations for the MMRT as a function of $n$
(100 simulations per time step $dt=0.01$ with $n=4$ to $15$ rods, $D=10$ for the
diffusion constant , $L=0.3$ for the distance from the origin
${\bf 0}$ and $l_{0}=0.25$ the length of each rod). The initial
configuration is such that $\left|
\varphi_{n}(0)
\right| <2\pi-\varepsilon$, e.g. the straight initial configuration:
$\theta_{k}(0)=0, \: \forall k=(1,...,n)$.
\begin{figure}
    \centering
    \includegraphics[width=0.60\textwidth]{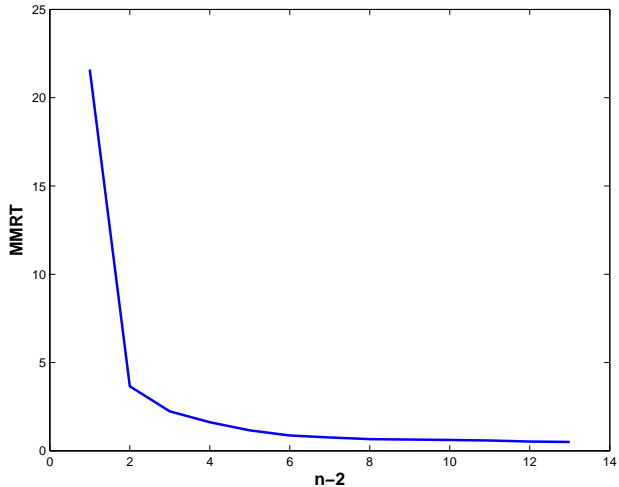}  \\
    \caption{ {\bf Brownian simulations of the MMRT as a function of the number of beads $n-2$ for
$D=10$.} }\label{figETvar}
\end{figure}
In Figure \ref{figDmanynL}, we present some Brownian simulations for
the MMRT as a function of $D$ and of $L$ (Figure
\ref{figDmanynL}(a) and Figure
\ref{figDmanynL}(b) respectively). $L$ and $l_{0}$ satisfy the rotation
compatibility condition $n l_{0}>L$. $E\left[\tau_{min}\right]$
decreases with $D$ and increases with the distance from the origin
$L$. It remains an open problem to compute the MMRT asymptotically for $n$ large.

\begin{figure}
\centering
\begin{tabular}{cc}
(a)\includegraphics[width=0.50\textwidth]{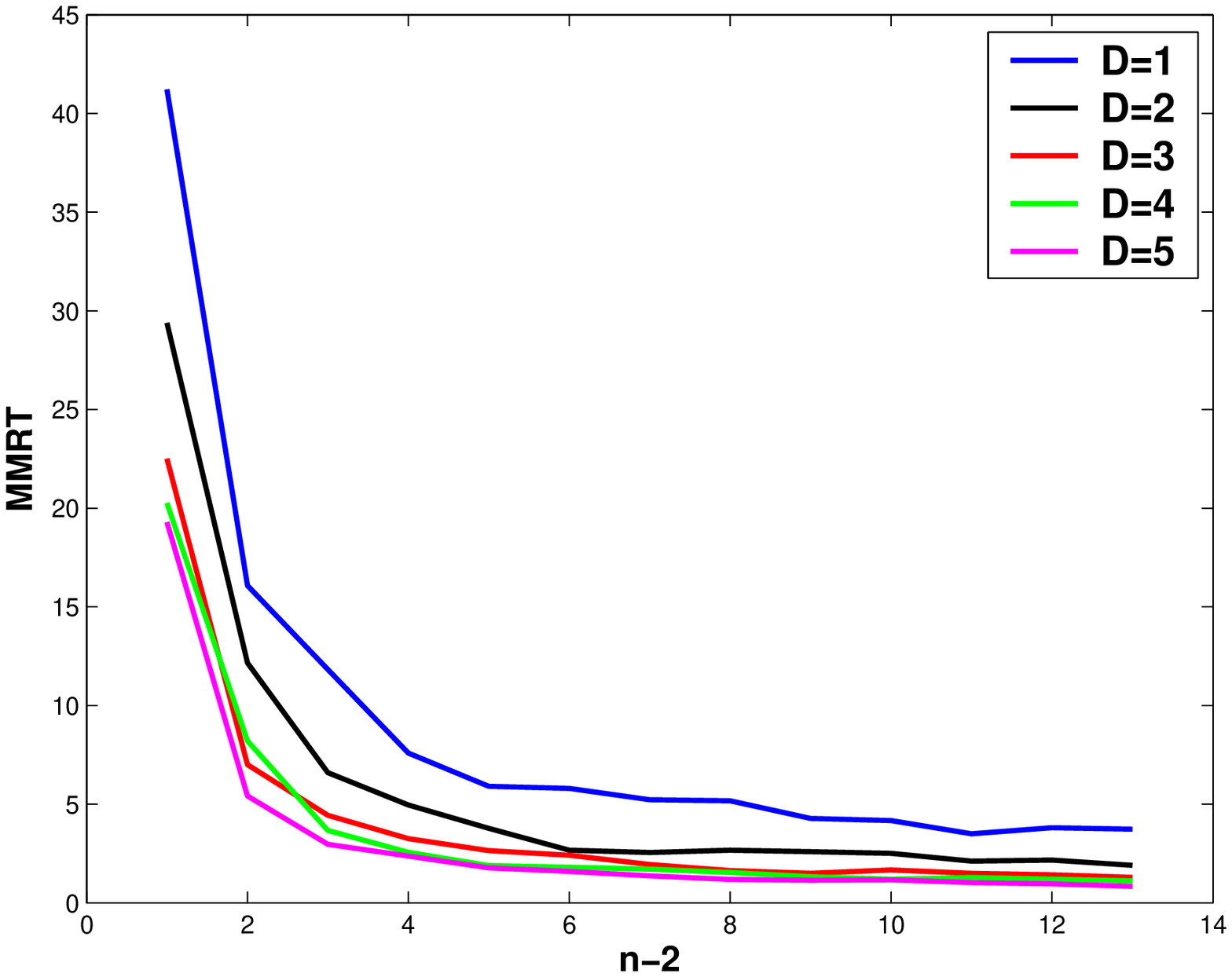} &
(b)\includegraphics[width=0.50\textwidth]{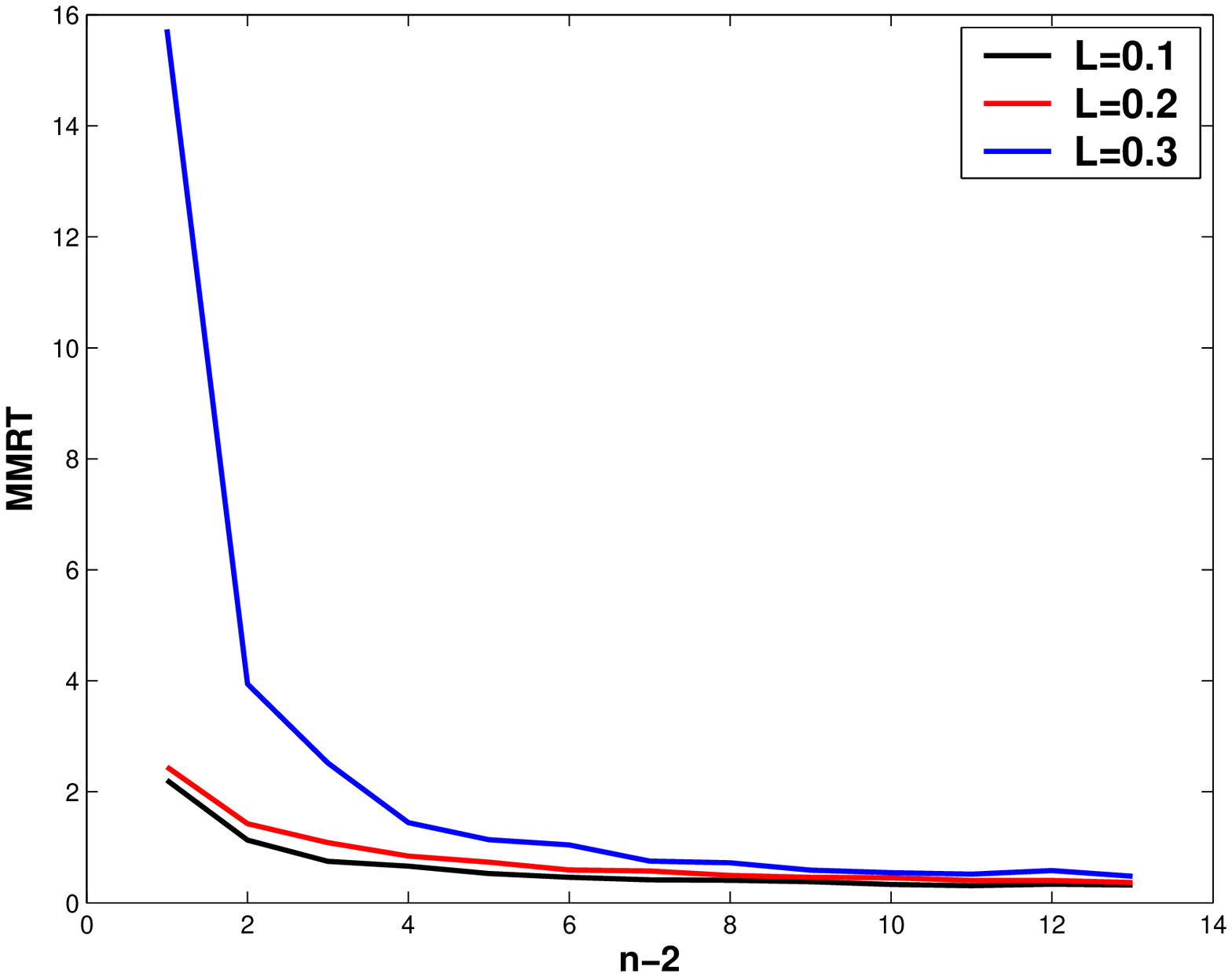} \\
\end{tabular}
\caption{ {\bf Brownian simulations of the
 MMRT as a function of the number of beads $n-2$. } (a)
 For several values of $D$ with $L=0.3$, (b) for several values of $L$ with $D=10$. }\label{figDmanynL}
\end{figure}

\subsection{Initial configuration in the boundary layer}\label{2inconf}
When the polymer has initially almost made a loop, we expect the MRT to have a different behavior. We look at this numerically, and we start with an
initial polymer configuration in the boundary layer :
\beqq
2\pi-\varepsilon \leq \left|
\varphi_{n}(0) \right| < 2\pi.
\eeqq
where $\varphi_n$ is defined in eq. (\ref{2XnRphi}). The rotation of the polymer will be completed very fast and using Brownian
simulations (100 runs per point with $n=10$ rods and $D=1$ for the diffusion constant, $L=0.1$ for the distance from the
origin ${\bf 0}$ and $l_{0}=0.2$ for the length of each rod), we plotted in Figure \ref{InitialTotalAngle} the results
showing  that when the initial total angle $\varphi_{n}(0)$ tends to $2\pi$, the MRT
tends to zero, with the precise asymptotic remaining to be completed. Numerically, we postulate that there is a threshold for an
initial total angle $\left| \varphi_{n}(0) \right|
=\frac{\pi}{2}$. When $\left| \varphi_{n}(0) \right|<\frac{\pi}{2}$,
the MRT appears to be independent from this angle, whereas for $\left| \varphi_{n}(0)
\right|>\frac{\pi}{2}$, the MRT decreases to zero. However, this needs further investigation.

\begin{figure}
\centering
\includegraphics[width=0.60\textwidth]{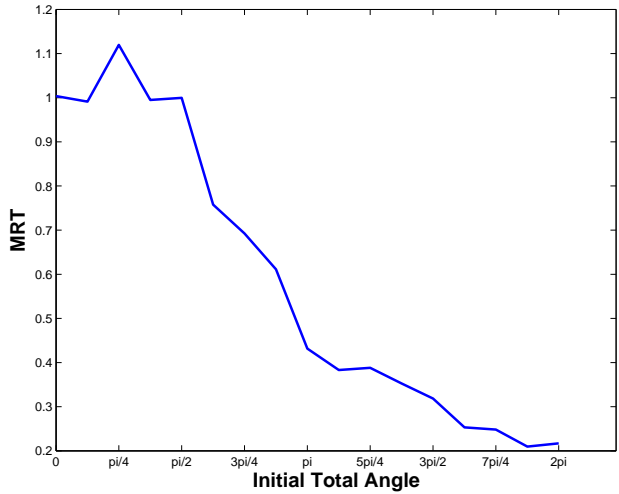}
\caption{ {\bf Brownian simulations of the
 MRT as a function of the initial total angle $\varphi_{n}(0)$, where $\varphi_{n}$ is defined in (\ref{2XnRphi}). } }\label{InitialTotalAngle}
\end{figure}

\section{Discussion and conclusion}
In the present paper, we studied the MRT for a planar random polymer
consisting of $n$ rods of length $l_0$. The first end is fixed at
a distance $L$ from the origin, while the other end moves as a
Brownian motion. Interestingly, we have shown here that the motion of the free polymer end satisfies a new
stochastic equation (\ref{2Xnfinal}), containing a nonlinear
time-dependent deterministic drift. When $n$ is large, the limit process is an Ornstein-Uhlenbeck process,
with different time scales for each of the two coordinates.

We found that the MRT $E\left[\tau_{n}\right]$ actually depends on the mean initial configuration of the polymer. This result is in contrast with the one of the small hole theory \cite{WaK93,WHK93,HoS04,SSHE06,SSH06a,SSH06b,SSH07,PRE2008} where the leading order term of the MFPT for a Brownian particle to reach a small hole does not depend on the initial configuration. Although the MRT is not falling exactly into the narrow escape problems, it is a rare event and the polymer completes a rotation when the free end reaches any point of the positive $x$-axis. The reason why the initial configuration survives in the large time asymptotic  regime, is due to the dynamics of the free moving end, approximated as a sum of i.i.d. variables, which is not Markovian. It leads to a process with memory. In summary, for $n
l_0>>L$, the leading order term of the MRT is given by:
\begin{enumerate}
    \item {for a general initial configuration:
\beqq
E\left[\tau_{n}\right] \approx  \frac{\sqrt{n}}{8D} \left[2\ln
\left(\frac{c_{n}}{\sqrt{n}}\right)+0.08\frac{n}{c_{n}^{2}}+Q\right],
\eeqq
where $c_{n}\equiv E \left( \sum^{n}_{k=1} e^{\frac{i}{\sqrt{2D}}
\theta_{k}(0)}\right)$,
$\left(\theta_{k}(0),1\leq k\leq n\right)$ is the sequence of the initial angles and $Q\approx9.54$,}
\item { for a stretched initial configuration:
\beqq
E\left[\tau_{n}\right] \approx  \frac{\sqrt{n}}{8D} \; \left(\ln (n)+Q\right),
\eeqq } \\
    \item {for an average over uniformly distributed initial angles:
\beqq
E\left[\tau_{n}\right] \approx  \frac{\sqrt{n}}{8D} \; \tilde{Q},
\eeqq
where $\tilde{Q}\approx9.62$. }
\end{enumerate}
As we have shown, these formulas are in very good agreement with Brownian
simulations (Fig. \ref{2comparisonimpr}). \\
After completion of this work, several questions arise naturally, namely:\\
- When the polymer has already made one loop, what is the probability that it makes a second loop before unwrapping and in that case, what is the MRT?  \\
- How can we extend our study in dimension 3 and higher dimensions? See \cite{HVY11} for some first results.

\appendix
\section{Appendix: Proof of Theorem \ref{CLT}} \label{apTCL}

We show the convergence of the sequence $(M_{t}^{(n)},t\geq 0)$ to a
Brownian motion in dimension 2, with a different time scale function for each coordinate.
This classically involves 2 steps (see e.g. \cite{Bil68,Bil78} or \cite{ReY99} (Chapter XIII)):
\begin{enumerate}
    \item {the convergence of the finite dimensional distributions, and}
    \item {the tightness of the sequence $(M_{t}^{(n)},t\geq0)$}
\end{enumerate}
$\left.a\right)$ The martingales $S_{t}^{(n)}$ and $C_{t}^{(n)}$ may be written in the form of stochastic integrals, as:
\beqq
S_{t}^{(n)} &=& \frac{1}{\sqrt{n}}  \int^{t}_{0} \sum^{n}_{k=1} \sin( B_{k}(s)) \: dB_{k}(s) \\
C_{t}^{(n)} &=& \frac{1}{\sqrt{n}} \int^{t}_{0} \sum^{n}_{k=1} \cos( B_{k}(s)) \: dB_{k}(s).
\eeqq
Consequently, with $<M>_{t}$ denoting the quadratic variation \cite{ReY99} of the martingale $\left(M_{t},t\geq0\right)$, we obtain:
\beq
    \left\langle S^{(n)}\right\rangle_{t} &=&  \frac{1}{n} \int^{t}_{0} \sum^{n}_{k=1} \sin^{2}( B_{k}(s)) \: ds
    \overset{{a.s.}}{\underset{n\rightarrow\infty}\longrightarrow} \int^{t}_{0} E\left[  \sin^{2}(B_{1}(s)) \right] \: ds, \label{2crochetS} \\
    \left\langle C^{(n)}\right\rangle_{t} &=&  \frac{1}{n} \int^{t}_{0} \sum^{n}_{k=1} \cos^{2}(B_{k}(s)) \: ds \overset{{a.s.}}{\underset{n\rightarrow\infty}\longrightarrow} \int^{t}_{0} E\left[  \cos^{2}(B_{1}(s)) \right] \: ds, \label{2crochetC} \\
    \left\langle S^{(n)},C^{(n)} \right\rangle_{t} &=&  -\frac{1}{2n} \int^{t}_{0} \sum^{n}_{k=1} \sin(2B_{k}(s)) \: ds
    \overset{{a.s.}}{\underset{n\rightarrow\infty}\longrightarrow} \int^{t}_{0} E\left[  \sin(2B_{1}(s)) \right] \: ds, \nonumber \\
    \label{2crochetSC}
\eeq
where the classical ergodic theorem yields the a.s. convergence \cite{Bil68}. We now give explicit expressions for the 3 right-hand sides. \\
Consider:
\beqq
E \left[\left(\exp(i B_{1}(s))\right)^{2}\right] = E \left[\exp(2i B_{1}(s))\right]=e^{-2s},
\eeqq
from which we deduce:
\beqq
E \left[ \cos^{2}( B_{1}(s)) - \sin^{2}( B_{1}(s)) \right]=E \left[ \cos( 2B_{1}(s))\right]=e^{-2s},
\eeqq
and
\beqq
E \left[ \sin( 2B_{1}(s))\right]=0.
\eeqq
Finally, we obtain:
\beq\label{2Ecossquared}
    E \left[ \cos^{2}( B_{1}(s)) \right] = \frac{1+e^{-2s}}{2},
\eeq
\beq\label{2Esinsquared}
    E \left[ \sin^{2}( B_{1}(s)) \right] = \frac{1-e^{-2s}}{2}.
\eeq
$\left.b\right)$ The previous results allow to obtain the convergence in law for $\left(M^{(n)}_{t_{1}},\ldots,M^{(n)}_{t_{k}}\right)$,
say, as $n\rightarrow\infty$.
Indeed, by taking $f,g:\mathbb{R}_{+}\rightarrow\mathbb{R}$, simple functions, we may write:
\beq
 E \left[ e^{i\left(\int^{\infty}_{0} f(u) dS^{(n)}_{u} + \int^{\infty}_{0} g(u) dC^{(n)}_{u} \right)} \right] \equiv E \left[ \exp\left(i \; \Sigma^{(n)}_{\infty}\right)\right],
\eeq
where:
\beq
 \Sigma^{(n)}_{t}= \int^{t}_{0} \left(f(u) dS^{(n)}_{u} + g(u) dC^{(n)}_{u} \right),
\eeq
as follows:
\beq\label{2characteristicfunction}
E \left[ \exp\left(i \; \Sigma^{(n)}_{\infty}\right)\right]=E \left[ \exp\left(i \; \Sigma^{(n)}_{\infty}+\frac{1}{2}\; <\Sigma^{(n)}>_{\infty}\right)\exp\left(-\frac{1}{2} \; <\Sigma^{(n)}>_{\infty}\right)\right].
\eeq
The results obtained in part $\left.a\right)$, together with the fact that \cite{ReY99}:
\beq
E \left[ \exp\left(i \; \Sigma^{(n)}_{\infty}+\frac{1}{2}\; <\Sigma^{(n)}>_{\infty}\right)\right]=1
\eeq
now yield:
\beq
E \left[ \exp\left(i \; \Sigma^{(n)}_{\infty}\right)\right]{\underset{n\rightarrow\infty}\longrightarrow} \exp\left( -\frac{1}{2} \int^{\infty}_{0} \left(f^{2}(u) \frac{1-e^{-2u}}{2} + g^{2}(u) \frac{1+e^{-2u}}{2}\right) \: du \right).
\eeq
$\left.c\right)$ It now remains to prove the tightness \cite{Bil78} of the distributions of the sequence $M^{(n)}$, which follows
from a classical application of Kolmogorov's criterion; indeed, for $\beta>0$ and $c_{\beta}$ a positive constant,
\beq\label{2tightness3}
    && E \left[ \left| M_{t}^{(n)}-M_{s}^{(n)} \right|^{2\beta} \right] \nonumber \\
    &\leq& c_{\beta} \left\{E \left[ \left(\left\langle S^{(n)}\right\rangle_{t}-\left\langle S^{(n)}\right\rangle_{s}\right)^{\beta} \right]
    + E \left[ \left(\left\langle C^{(n)}\right\rangle_{t}-\left\langle C^{(n)}\right\rangle_{s}\right)^{\beta} \right]\right\} \nonumber \\
    &\leq& 2 c_{\beta} \left| t-s \right|^{\beta}.
\eeq
We refer the reader who may want more details about the arguments used to \cite{ReY99}, Chapter XIII, where convergence
in distribution on the canonical space $C(\mathbb{R}_{+},\mathbb{R})$ is discussed. \hfill \QED

\newpage


\end{document}